\documentclass[10pt,  twoside, onecolumn, notitlepage]{article}
\usepackage[vmarginratio=1:1,a4paper,body={6.5in, 9.5in}]{geometry}
\usepackage[utf8x]{inputenc}

\usepackage{amsmath,amssymb}
\usepackage{graphicx}

\title{Illustration of iterative linear solver behavior on simple $1D$ and $2D$ problems}
\author{Nicolas Ray\\
	INRIA \\
	\and
	Dmitry Sokolov \\
	Université de Lorraine \\
	}

\date{\today}
\begin{document}
\maketitle

\begin{abstract}
In geometry processing, numerical optimization methods often involve solving sparse linear systems of equations. These linear systems have a structure that strongly resembles to adjacency graphs of the underlying mesh. We observe how classic linear solvers behave on this specific type of problems. For the sake of simplicity, we minimise either the squared gradient or the squared Laplacian, evaluated by finite differences on a regular $1D$ or $2D$ grid. 

We observed the evolution of the solution for both energies, in $1D$ and $2D$, and with different solvers: Jacobi, Gauss-Seidel, SSOR (Symmetric successive over-relaxation) and CG (conjugate gradient \cite{Shewchuk}). Plotting results at different iterations allows to have an intuition of the behavior of these classic solvers.
\end{abstract}

\section*{Introduction}

When facing an optimization problem, the simplest approach is often to iterate, locally improving current solution. For example, heat diffusion can be done by iteratively replacing each value by the average value of its neighbors. With more background on numerical optimization, one would like to formulate this optimization as a linear problem. With such abstraction, it is possible to use solvers as black boxes: we do not have an interpretation of the evolution of the solution during the iterations.

We are interested in observing the behavior of basic linear solvers for very simple cases. It does not necessarily impact the way to use them, but provides some intuition about the magic of linear solvers in the context of minimizing energies on a mesh. 

Our meshes are regular $1D$ and $2D$ grids of size $n$, and variables are attached to vertices. We address the coordinates of the unknown vector $\mathbf{x}$ as follows: in $1D$ $\mathbf{x}_i$ is the $i^{th}$ element of $\mathbf{x}$ and in $2D$, the value $\mathbf{x}_{i,j}$ of the vertex located at the $i^{th}$ line and $k^{th}$ column is the $(n*i+j)^{th}$ element of $\mathbf{x}$. 

Our tests are performed on the minimization of two energies:
\begin{itemize}
\item {\it Gradient energy:} The sum of squared difference between adjacent samples: $\Sigma_i(\mathbf{x}_i-\mathbf{x}_{i+1})^2$ in $1D$ and $\Sigma_{i,j}(\mathbf{x}_{i,j}-\mathbf{x}_{i+1,j})^2 +\Sigma_{i,j} (\mathbf{x}_{i,j}-\mathbf{x}_{i,j+1})^2$ in $2D$. 
\item {\it Laplacian energy:} The sum of squared difference between a sample and the average of its neighbors: $\Sigma_i (\mathbf{x}_i-\mathbf{x}_{i-1}/2-\mathbf{x}_{i+1}/2)^2$ in $1D$. In $2D$, it is the same on energy for boundary vertices, and $\Sigma_{i,j} (\mathbf{x}_{i,j}-\mathbf{x}_{i-1,j}/2-\mathbf{x}_{i+1,j}/2)^2 + (\mathbf{x}_{i,j}-\mathbf{x}_{i,j-1}/2-\mathbf{x}_{i,j+1}/2)^2$ over other vertices.
\end{itemize}

Note that for the Laplacian energy with the right number of constraints ($2$ in $1D$, and $3$ in $2D$), the minimization of the first energy is equivalent to the minimization of the second energy.

\section{$1D$ problems}

Observing $1D$ problems allows to visualize the convergence by a $3D$ surface. In all figures in this section, the first axis (horizontal) represents the domain --- a $1D$ grid, the second axis (vertical) represents the value associated to each vertex, and the third axis (depth) is the time.

\subsection{Gradient energy}

The $1D$ grid has $n$ vertices, here we develop the matrices for $n=5$. The problem is formalized by this set of equations $Ax=b$~:
$$
\begin{pmatrix}
   0 & -1 & 0 & 0 & 0 \\
   0 & 1 & -1 & 0 & 0 \\
   0 & 0 & 1 & -1 & 0 \\
   0 & 0 & 0 & 1 & 0 \\
   1 & 0 & 0 & 0 & 0 \\
   0 & 0 & 0 & 0 & 1 
\end{pmatrix}
\mathbf{x}= 
\begin{pmatrix}
   -2\\
   0 \\
   0 \\
   6\\
   2 \\
   6
\end{pmatrix}
$$

The $4$ first lines represent the objective. To enforce boundary constraints $\mathbf{x}_1=2$ and $\mathbf{x}_5=6$, we replaced  coefficients applied to $\mathbf{x}_1$ and $\mathbf{x}_5$ by a contribution to the right hand side. We also add the last two lines to explicitly set $\mathbf{x}_1=2$ and $\mathbf{x}_5=6$. In the least squares sense, it produces the new system to solve $A^\top A\mathbf{x}=A^\top b$ i.e.

$$
\begin{pmatrix}
   1 & 0 & 0 & 0 & 0 \\
   0 & 2 & -1 & 0 & 0 \\
   0 & -1 & 2 & -1 & 0 \\
   0 & 0 & -1 & 2 & 0 \\
   0 & 0 & 0 & 0 & 1
\end{pmatrix}
\mathbf{x}= 
\begin{pmatrix}
   2\\
   2 \\
   0 \\
   6\\
   6
\end{pmatrix}
$$

The condition number of $A^\top A$ is $5.9$ for $n=5$, and $1000$ for $n=50$. These numbers are not high enough to produce numerical instabilities.

\begin{figure}[h]
\centerline{\includegraphics[width = 16cm]{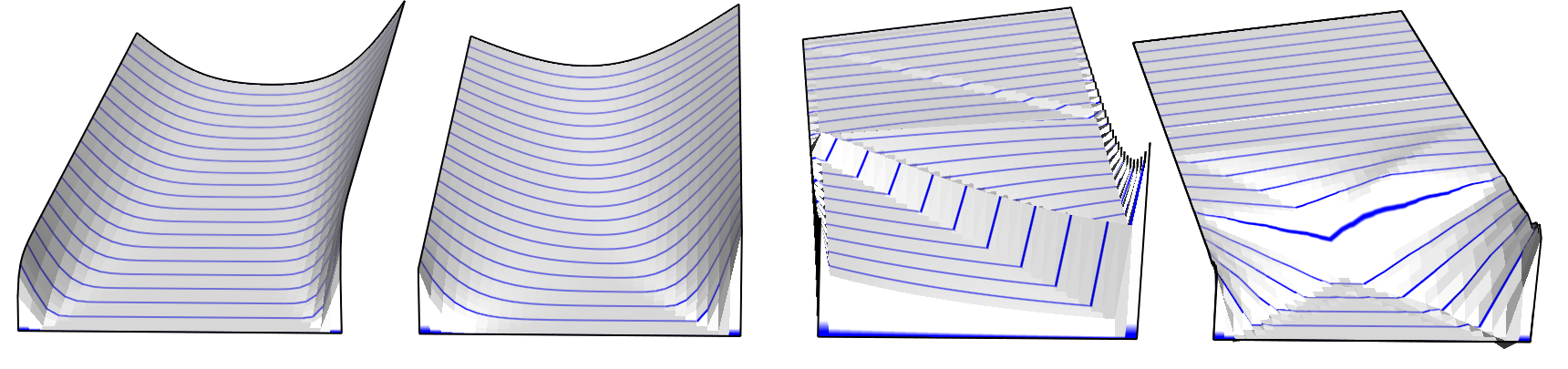}}
\centerline{Both extremities are locked to strictly positive values.}
\centerline{\includegraphics[width = 16cm]{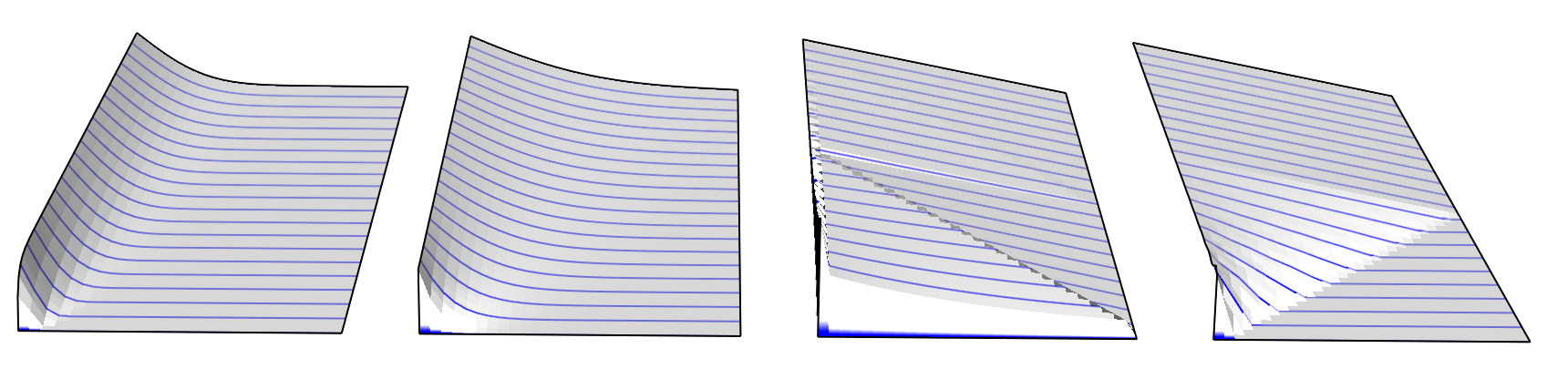}}
\centerline{The left extremity is locked to a strictly positive value, right to zero.} 
\centerline{\includegraphics[width = 16cm]{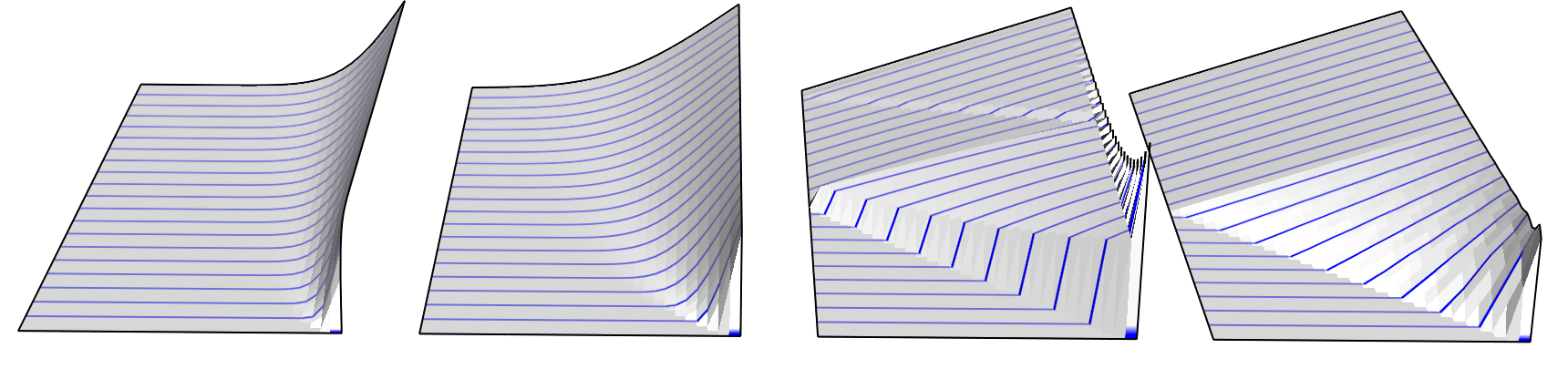}}
\centerline{The right extremity is locked to a strictly positive value, left to zero.}
\caption{Evolution of the solution for the $1D$ gradient energy with $20$ variables solved by respectively Jacobi, Gauss-Seidel, SSOR and CG with 120 iterations. The horizontal dimension is the $1D$ grid of size $20$, the vertical dimension is the value of our variables, and the depth dimension is the time. Each blue curve corresponds to the solution at a given iteration.}
   \label{fig:1Dconvergence}
\end{figure}

This first experience (Fig.~\ref{fig:1Dconvergence}) shows the behavior of each solver with the same matrix. The boundary conditions (locked coordinates of $\mathbf{x}$) are set in the way that the solution is the straight line defined by $\mathbf{x}_i = i(\mathbf{x}_n-\mathbf{x}_1)/n$. We observe the expected relative speed of convergence i.e. Jacobi $<$ Gauss-Seidel $<$ SSOR $<$ Conjugate Gradient. We also visualize the impact of the order of the coordinates in Gauss-Seidel and SSOR: when only the left side is constrained (locked $\mathbf{x}_1$), the first iteration propagates the constraints to all coordinates whereas constraining only $\mathbf{x}_n$ take $n$ iteration to affect all the variables.

\begin{figure}[h]
\centerline{\includegraphics[width = 8cm]{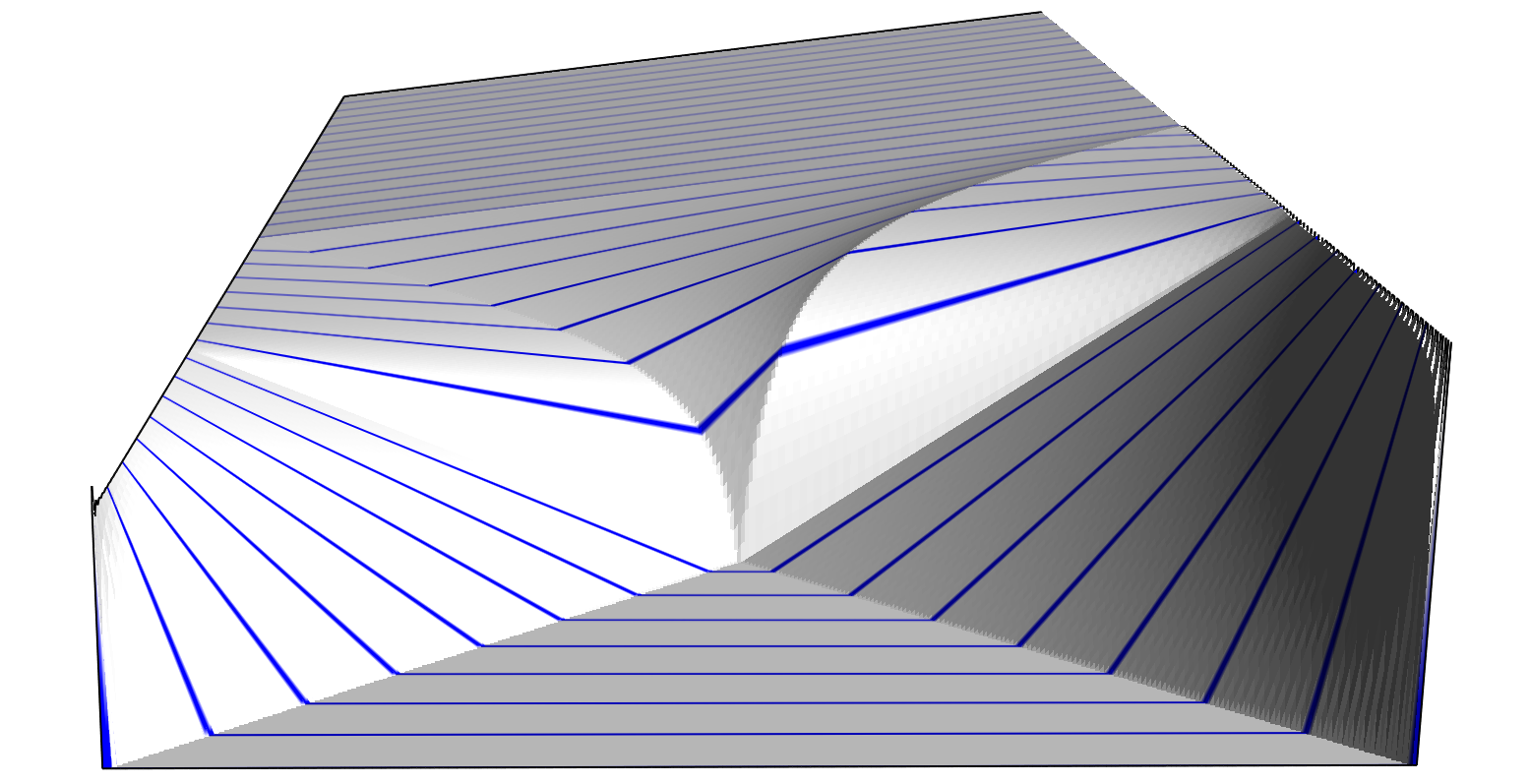}
\includegraphics[width = 8cm]{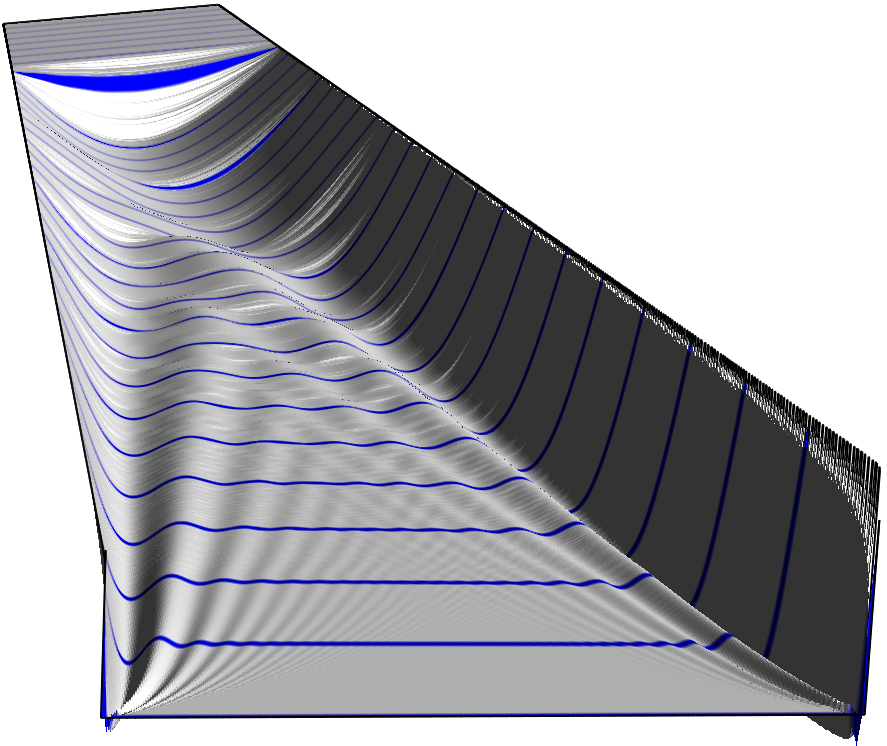}
}
\caption{Evolution of the solution for the $1D$ gradient energy (left) and the $1D$ laplacian energy (right) with the CG solver. The unknowned vector $\mathbf{x}$ is of dimension $100$ and the solvers did $500$ iterations. As usual, the horizontal dimension is the $1D$ domain, the vertical dimension is the value of our variables, and the depth dimension is the time.}
   \label{fig:CG1D}
\end{figure}

\subsection{Laplacian energy}

As for the gradient energy case, the $1D$ grid has $n$ vertices, and we develop the matrices for $n=5$. The problem is formalized by this set of equations $A\mathbf{x}=b$~:
$$
\begin{pmatrix}
   0 & 2 & -1 & 0 & 0 \\
   0 & -1 & 2 & -1 & 0 \\
   0 & 0 & -1 & 2 & 0 \\
   1 & 0 & 0 & 0 & 0 \\
   0 & 0 & 0 & 0 & 1 
\end{pmatrix}
\mathbf{x}= 
\begin{pmatrix}
   -2\\
   0 \\
   6\\
   2 \\
   6
\end{pmatrix}
$$

The $3$ first lines represent the objective. As for the ``gradient energy'', we replaced  coefficients applied to $\mathbf{x}_1$ and $\mathbf{x}_5$ by a contribution to the right hand side to enforce $\mathbf{x}_1=2$ and $\mathbf{x}_5=6$. We also add the last two lines to explicitly set $\mathbf{x}_1=2$ and $\mathbf{x}_5=6$. In least squares sense, it produces the new system to solve $A^\top A\mathbf{x}=A^\top b$ i.e.

$$
\begin{pmatrix}
   1 & 0 & 0 & 0 & 0 \\
   0 & 5 & -4 & 1 & 0 \\
   0 & -4 & 6 & -4 & 0 \\
   0 & 1 & -4 & 5 & 0 \\
   0 & 0 & 0 & 0 & 1
\end{pmatrix}
\mathbf{x}= 
\begin{pmatrix}
   2\\
   4 \\
   -8 \\
   12\\
   6
\end{pmatrix}
$$

The condition number of $A^\top A$ is $34.4$ for $n=5$, and $920000$ for $n=50$ coordinates. Such numbers produce numerical instabilities slowing down the solvers.

\begin{figure}[h]
\centerline{iteration $40$ \includegraphics[width = 14cm]{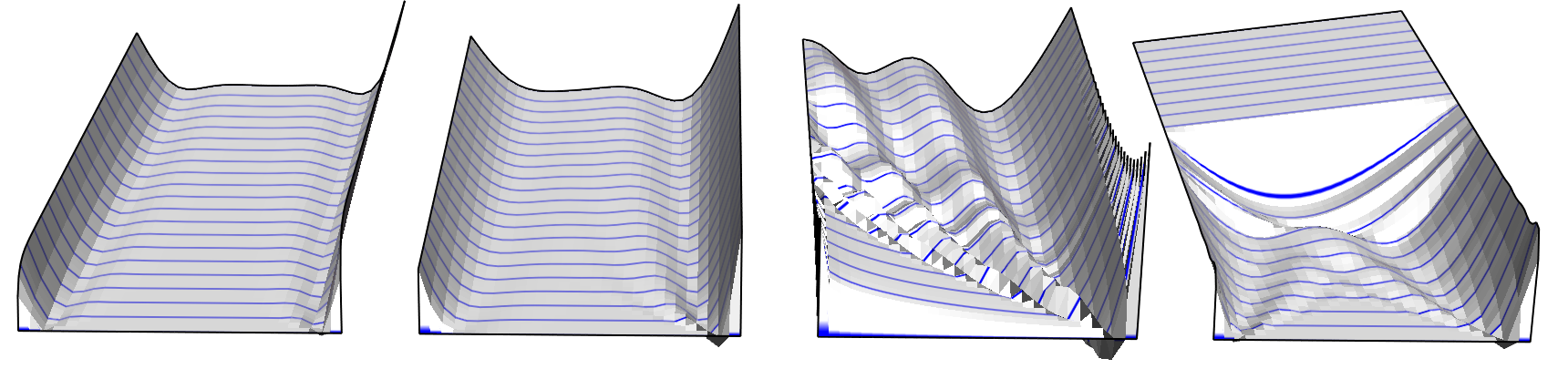}}
\centerline{iteration $4000$ \includegraphics[width = 14cm]{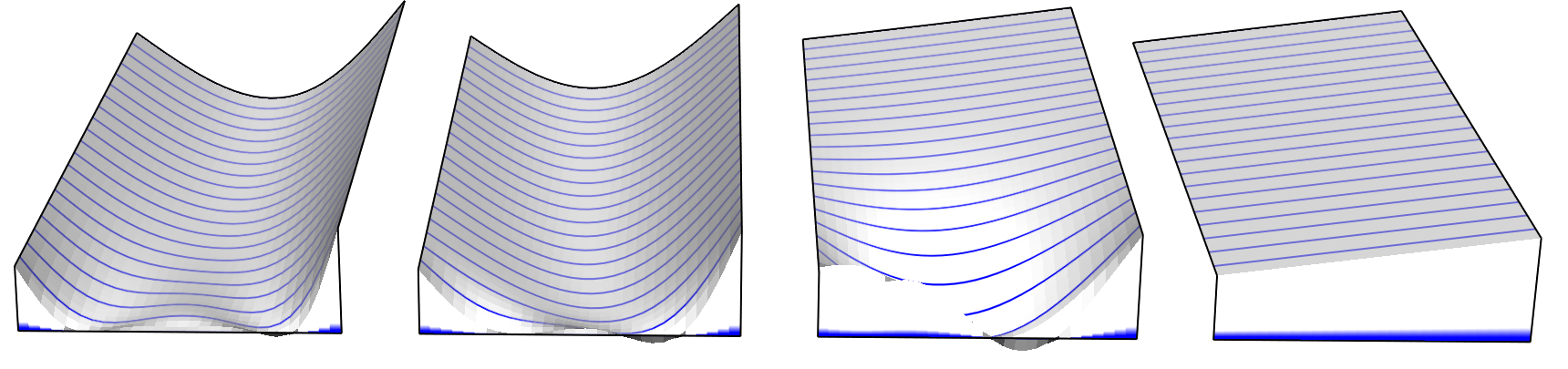}}
\centerline{iteration $40000$ \includegraphics[width = 14cm]{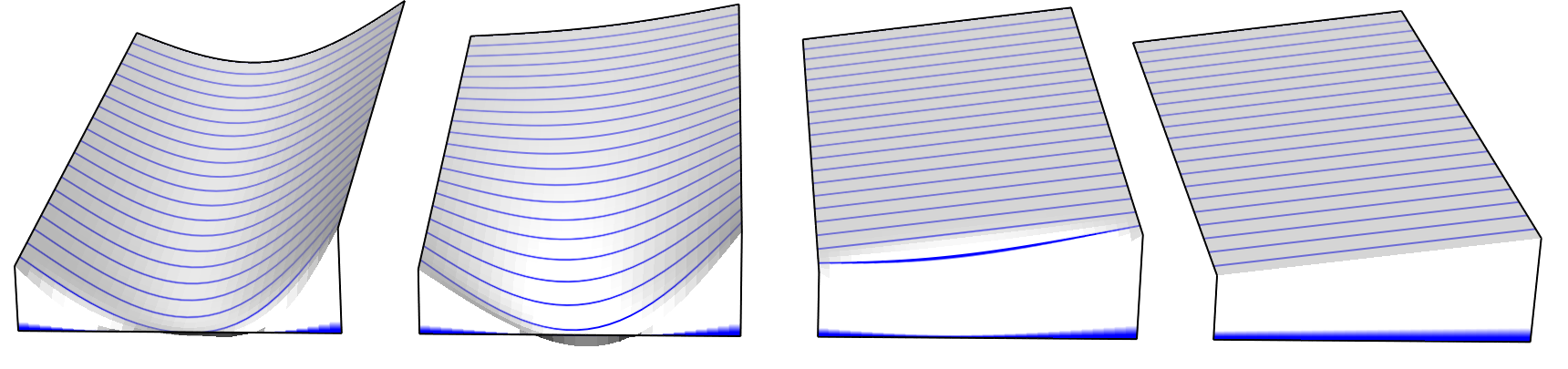}}
\centerline{--------------------- Jacobi  ------------------------ Gauss-Seidel ------------- SSOR ---------- Conjugate Gradient --}
\caption{Evolution of the solution for the $1D$ laplacian energy with $20$ variables solved by respectively Jacobi, Gauss-Seidel, SSOR and CG with $40$, $4000$, and $40000$ iterations. The horizontal dimension is the $1D$ grid of size $20$, the vertical dimension is the value of our variables, and the depth dimension is the time.}
   \label{fig:laplaceconvergence}
\end{figure}

In Fig.~\ref{fig:laplaceconvergence}, we observe the same expected relative speed of convergence i.e. Jacobi $<$ Gauss-Seidel $<$ SSOR $<$ Conjugate Gradient. Oscillations of the solution with Gauss-Seidel and SSOR are interesting to observe: it converges faster with higher local oscillations (in the solution space, not in the time dimension). 

Another interesting observation is that CG required more than $20$ iterations to converge. This was unexpected because, while CG is mostly used as an iterative solver, it is also a direct solver that is supposed to converge in a worse $n$ iterations. Our interpretation is that the poor conditionning of the matrix generates numerical instabilities that slow down the solver. This phenomena is even worst with higher dimension problem ($n=100$) as illustrated in Fig.~\ref{fig:CG1D}--right.

\section{$2D$ problems}

\begin{figure}[h]
\centerline{Iteration 1 \includegraphics[width = 14cm]{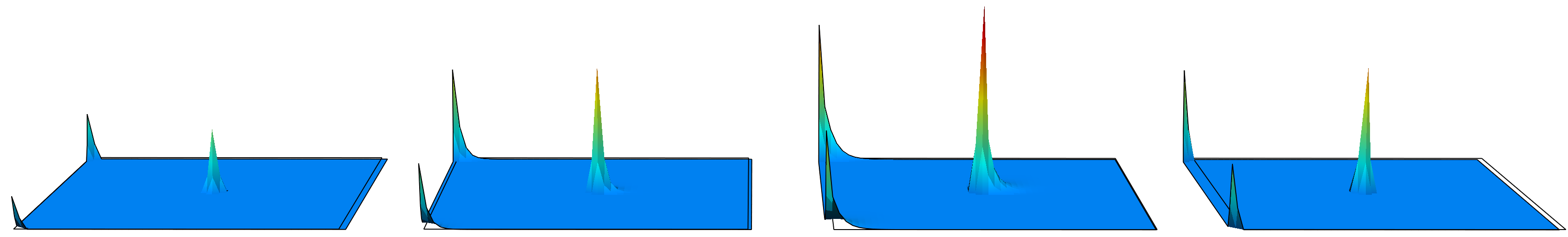}}
\centerline{Iteration 5 \includegraphics[width = 14cm]{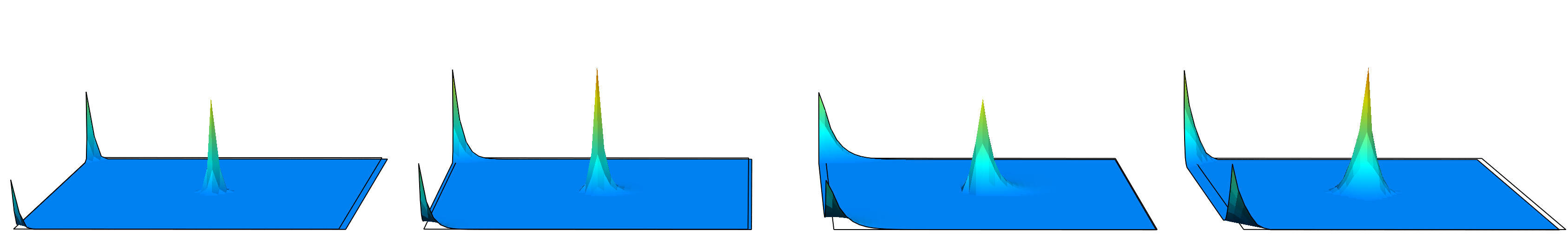}}
\centerline{Iteration 10 \includegraphics[width = 14cm]{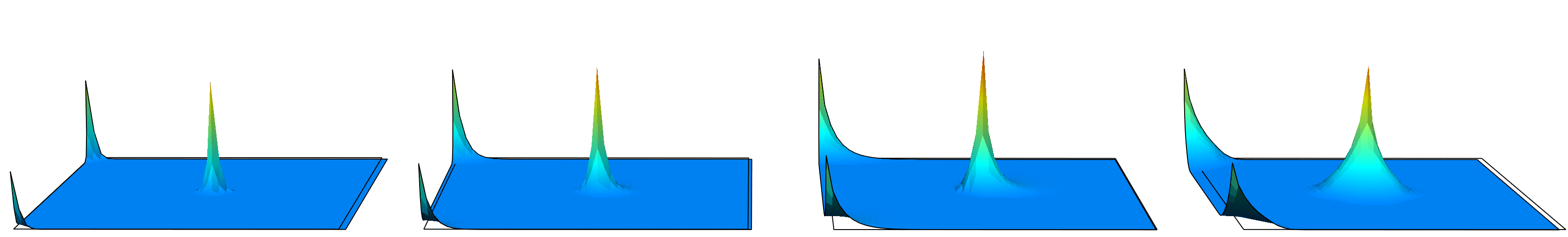}}
\centerline{Iteration 50 \includegraphics[width = 14cm]{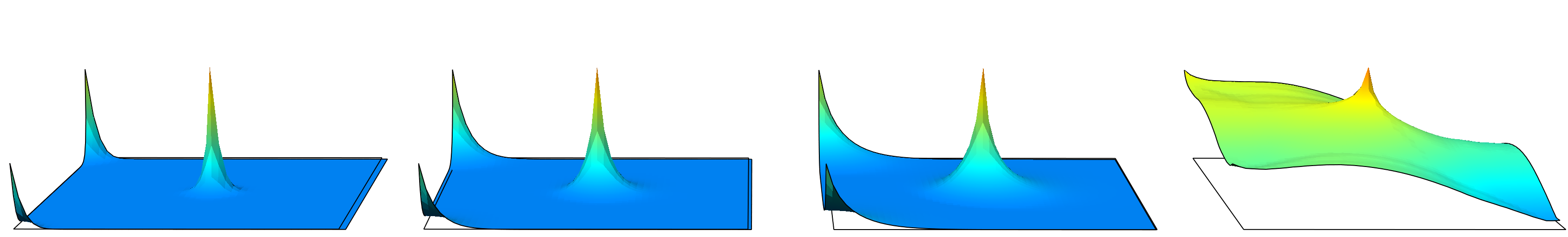}}
\centerline{Iteration 100 \includegraphics[width = 14cm]{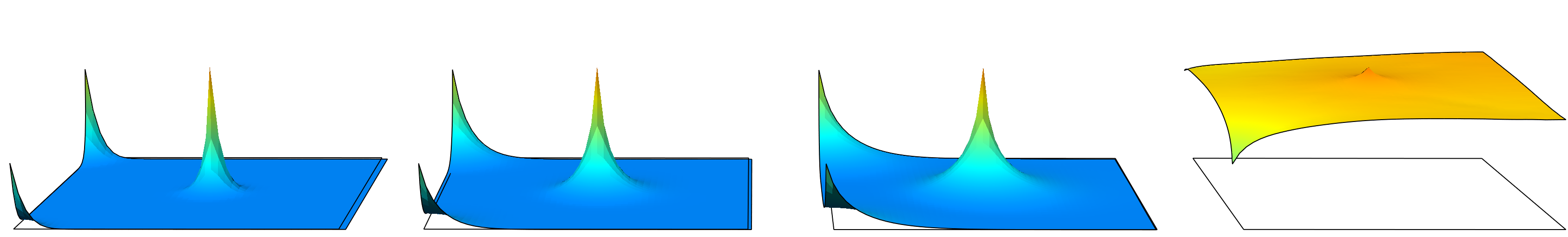}}
\centerline{Iteration 2400 \includegraphics[width = 14cm]{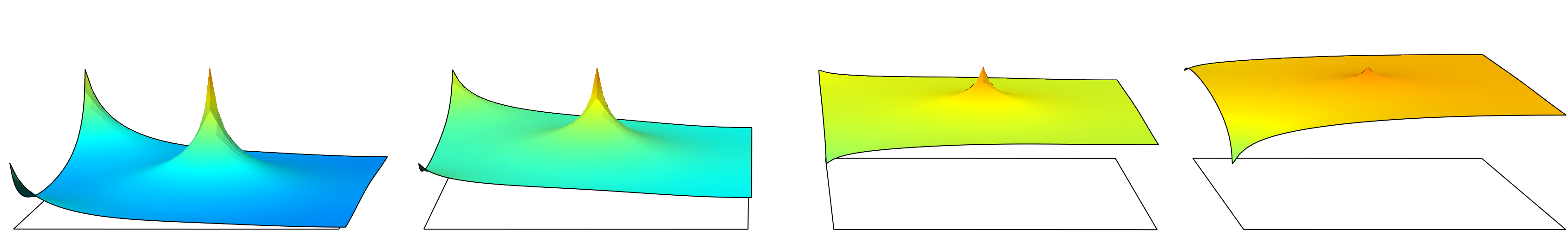}}
\centerline{Iteration 2500 \includegraphics[width = 14cm]{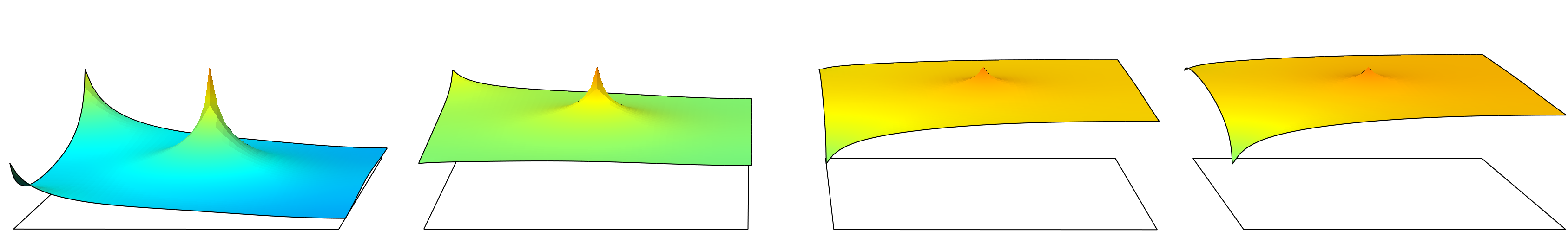}}
\centerline{Iteration 10000 \includegraphics[width = 14cm]{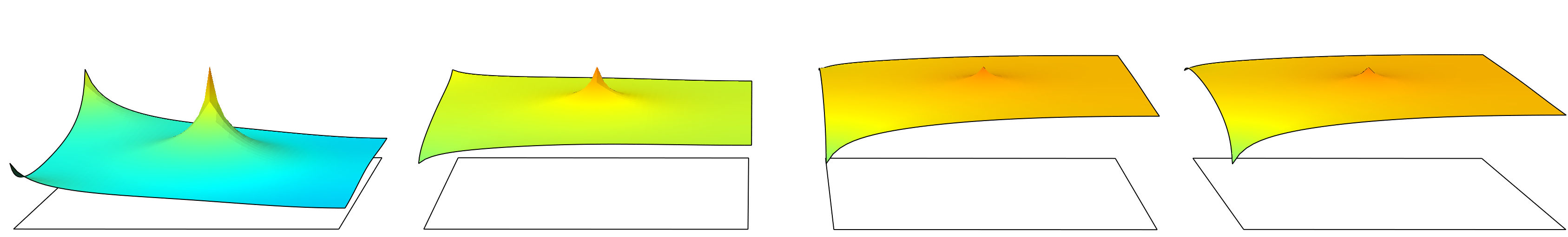}}
\centerline{------------------------------ Jacobi  -------- Gauss-Seidel -------- SSOR ------ Conjugate Gradient --}
\caption{Plots of the solution of the $2D$ problem with the gradient energy}
   \label{fig:2Dgradient}
\end{figure}

\begin{figure}[h]
\centerline{Iteration 1 \includegraphics[width = 14cm]{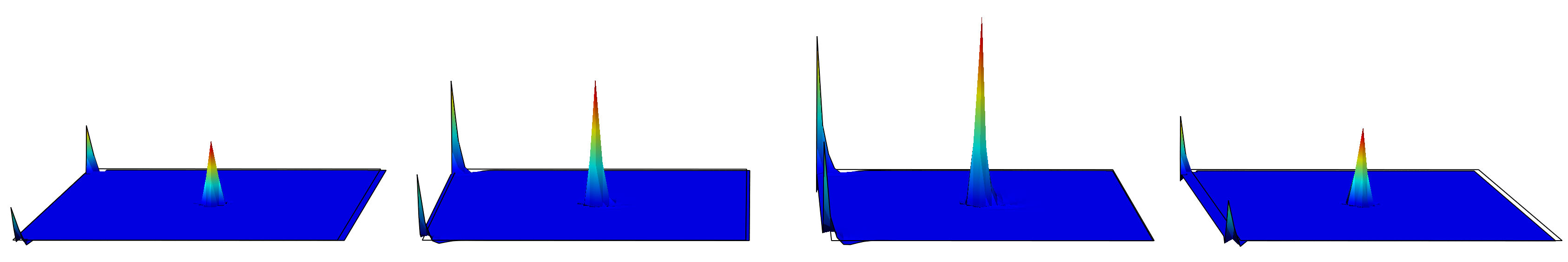}}
\centerline{Iteration 5 \includegraphics[width = 14cm]{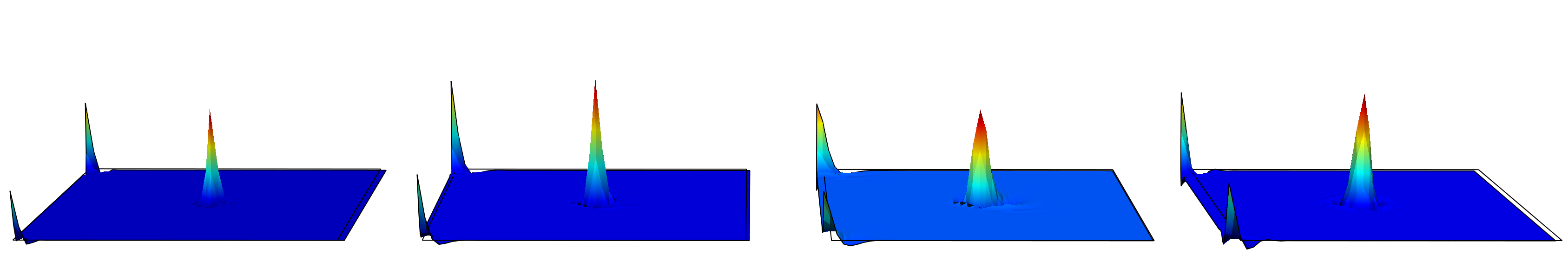}}
\centerline{Iteration 10 \includegraphics[width = 14cm]{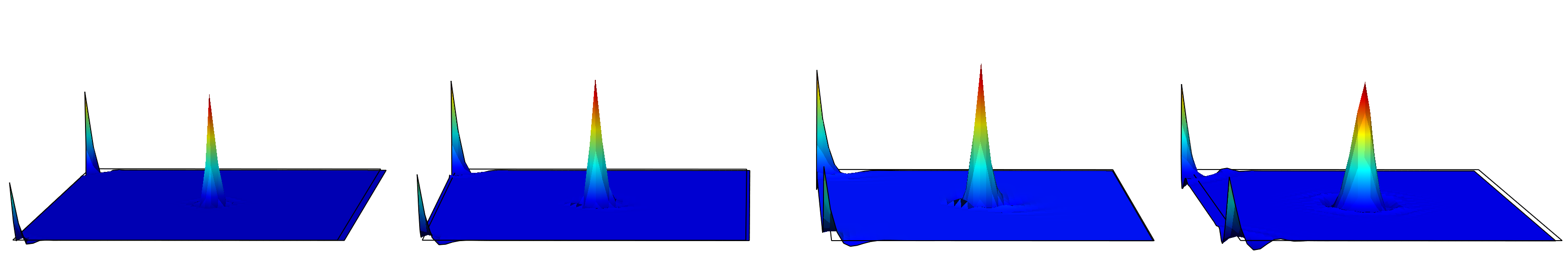}}
\centerline{Iteration 50 \includegraphics[width = 14cm]{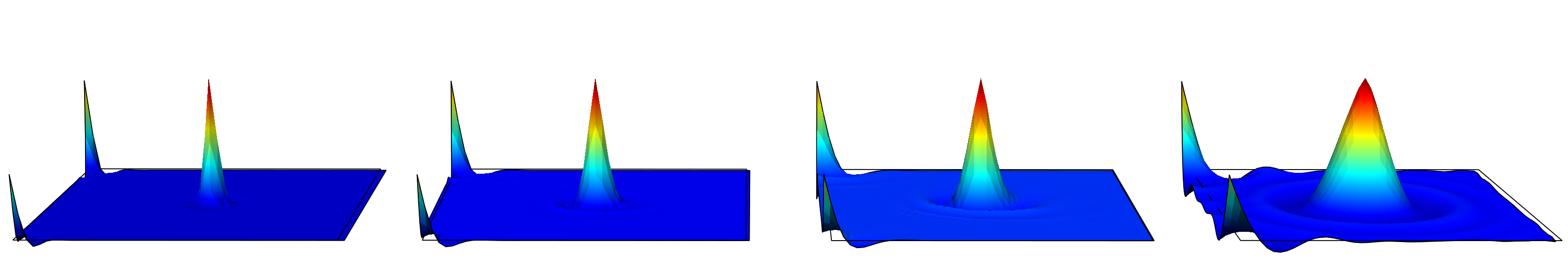}}
\centerline{Iteration 100 \includegraphics[width = 14cm]{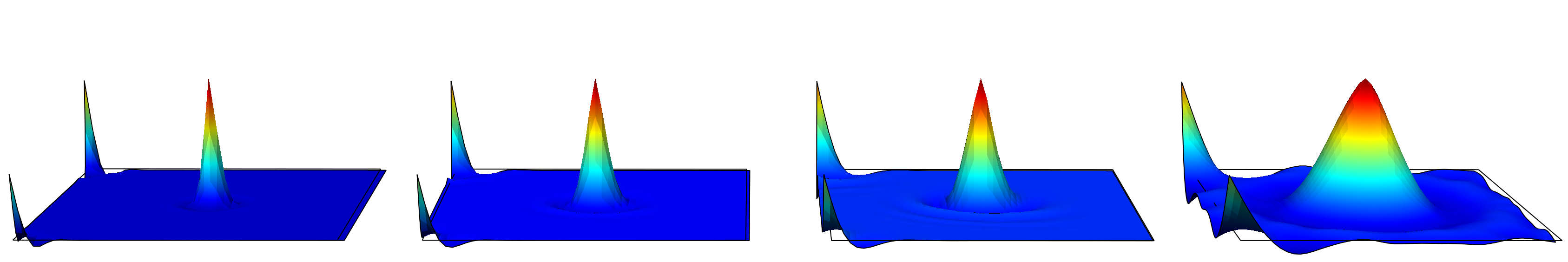}}
\centerline{Iteration 2400 \includegraphics[width = 14cm]{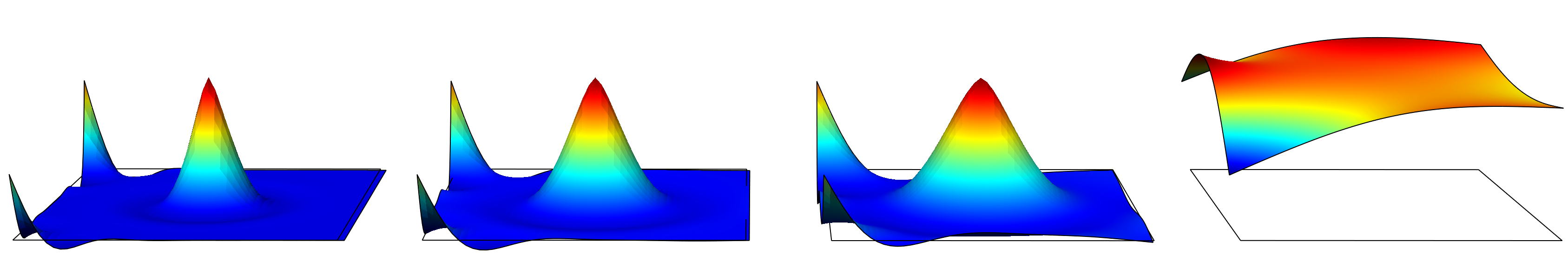}}
\centerline{Iteration 2500 \includegraphics[width = 14cm]{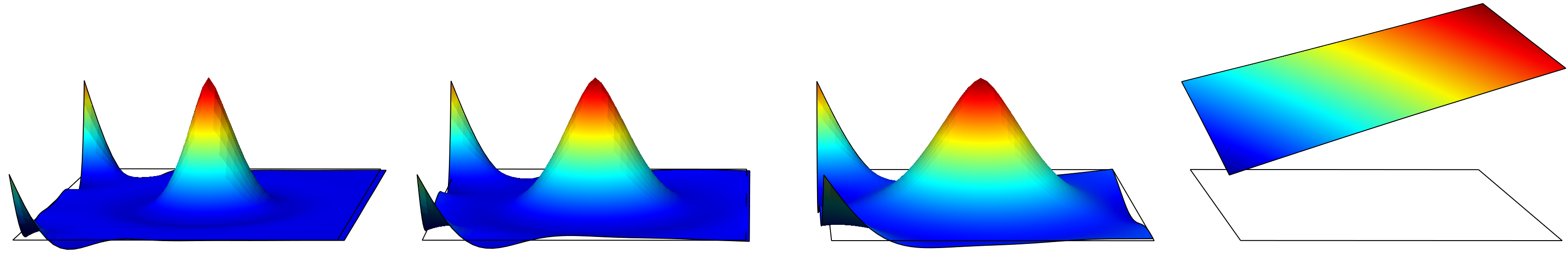}}
\centerline{Iteration 10000 \includegraphics[width = 14cm]{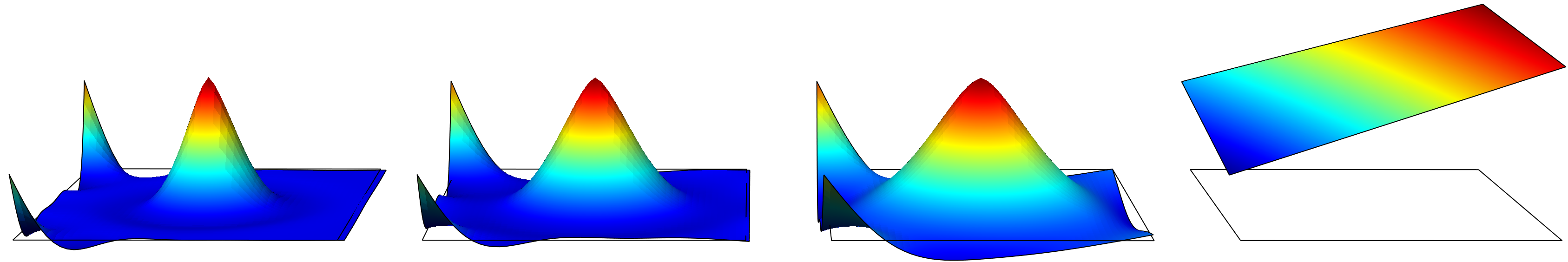}}
\centerline{------------------------------ Jacobi  -------- Gauss-Seidel -------- SSOR ------ Conjugate Gradient --}
\caption{Plots of the solution of the $2D$ problem with the Laplacian energy}
   \label{fig:2DLaplacien}
\end{figure}

The $2D$ problems are very similar to the $1D$ problems: we set the linear system $A\mathbf{x}=b$ with the gradient (resp. Laplacian) equations and add $3$ constraints ($2$ at corners of the grid and in the center). To solve it in the least squares senses, we have to solve the linear system $A^\top A\mathbf{x}=A^\top b$. We obtained results presented in Fig.~\ref{fig:2Dgradient} and \ref{fig:2DLaplacien}.

The observations are also very similar to $1D$: the speed of convergence is as usual i.e. Jacobi $<$ Gauss-Seidel $<$ SSOR $<$ CG as expected, and the non zero coordinates of $\mathbf{x}$ are "discovered" at each iteration in a front propagation fashion. Each constraint has an impact on the value of all vertices after $100$ iterations ($L^1$ distance of the diagonal). For the gradient energy, the CG method already have a fair solution and could be stopped at step $100$. For the Laplacian energy, that has a bad conditioning, the result at step $100$ is still far from the solution. Even at step $2400$, it is not converged... but at step $2500$, it's done ! (as expected).

We can also observe that the conditioning affects all methods (not only CG).


%
%

\section*{Conclusion: what did I learn?}

We observed two non obvious behaviors of the considered iterative linear solver:
\begin{itemize}
\item When starting with a null solution vector $\mathbf{x}$, all tested iterative solvers (including CG) remove the zero coefficient by front propagation from the constraints. It can be explained by the position of non zero coefficients in the matrix $A$ that is similar to the one of the vertices adjacency matrix. 
\item Ill-conditioned systems may take more than $n$ iterations to converge with the conjugate gradient algorithm. Stopping them earlier does not always give a good approximation of the solution e.g. $1D$ Laplacian.
\end{itemize}

\bibliographystyle{alpha}
\bibliography{CGbehavior}

\begin{thebibliography}{She94}

\bibitem[She94]{Shewchuk}
Jonathan~Richard Shewchuk.
\newblock An introduction to the conjugate gradient method without the
  agonizing pain,, 1994.
\newblock http://www.cs.cmu.edu/~quake-papers/painless-conjugate-gradient.pdf.

\end{thebibliography}
\end{document}